\documentclass[preprint]{imsart}

\usepackage{color, geometry}
\usepackage{amssymb}
\usepackage{hyperref}
\usepackage{amsthm}
\usepackage{cite}
\usepackage{float}
\usepackage{ulem} % REMOVE for submission!

\usepackage{graphicx}

\usepackage{fullpage}

\newcommand{\R}{\mathbb{R}}

\begin{document}

%%%%%%%%%%%%%%%%%% title page information %%%%%%%%%%%%%%%%%%
\begin{frontmatter}

\title{Expectation Maximization and the retrieval of the atmospheric extinction coefficients by inversion of Raman lidar data}

\begin{aug}
\author{Sara Garbarino,$^1$ Alberto Sorrentino,$^{2,3}$ Anna Maria Massone,$^3$ Alessia Sannino,$^4$ Antonella Boselli,$^5$ Xuan Wang,$^6$ Nicola Spinelli$^4$ and Michele Piana$^{2,3}$}

\address{$^1$Centre for Medical Image Computing, Department of Computer Science, University College London, UK\\
$^2$Dipartimento di Matematica, Universit\`a di Genova, Italy \\
$^3$CNR--SPIN, Genova, Italy \\
$^4$Dipartimento di Fisica, Universit\`a di Napoli Federico II, Italy\\
$^5$CNR--IMAA, Potenza, Italy \\
$^6$CNR--SPIN, Napoli, Italy
}
\end{aug}


\begin{abstract}
We consider the problem of retrieving the aerosol extinction coefficient from Raman lidar measurements. 
This is an ill--posed inverse problem that needs regularization, and we propose to use the
Expectation--Maximization (EM) algorithm to provide stable solutions. Indeed, EM is an iterative algorithm that
imposes a positivity constraint on the solution, and provides regularization if iterations are stopped early enough.
We describe the algorithm and propose a stopping criterion inspired by a statistical principle. 
We then discuss its properties concerning the spatial resolution.
Finally, we validate the proposed approach by using both synthetic data and experimental measurements; 
we compare the reconstructions obtained by EM with those obtained by the Tikhonov method, by the Levenberg-Marquardt method, 
as well as those obtained by combining data smoothing and numerical derivation.
\end{abstract}
\end{frontmatter}

%%%%%%%%%%%%%%%%%%%%%%%%%%  body  %%%%%%%%%%%%%%%%%%%%%%%%%%

\section{Introduction}\label{sec:sec1}

Lidar systems are increasingly used to characterize the temporal and spatial distribution of the aerosol optical characteristics in the atmosphere 
\cite{karyampudi1999validation,di2012raman}. Further, since lidar data present high spatial and temporal resolution, air motion can be precisely monitored with them \cite{kolev1988lidar}. However, lidar signals provide a very indirect measure of the relative concentration and distribution of the aerosols: the energy observed by a lidar is a function of the extinction and backscattering coefficients  which in turn are functions of the microphysical properties of the aerosols. Deriving such properties is therefore possible, but only via the numerical solution of two inverse problems: the first one allows the estimation of the extinction coefficients from the experimental data \cite{ansmann1990measurement,pappalardo2004aerosol,pornsawad2012retrieval,povey2014retrieval} and the second one allows the reconstruction of the size distribution, shape and refractive index of the aerosols \cite{muller1999microphysical,wang2007retrieval,osterloh2009parallel,osterloh2013regularized}, starting from the knowledge of the extinction and backscatter coefficient at multiple wavelengths.

The present paper introduces the use of Expectation Maximization (EM) for the solution of the first lidar inverse problem and shows that this iterative scheme is highly competitive, in terms of reconstruction accuracy, for this kind of data analysis. EM provides the unique non-negative solution of the maximum likelihood problem when the model equation is linear and the distance between the experimental and reproduced data is measured by means of a specific topology named the Kullback-Leibler (KL) distance \cite{bertero1998introduction}. The lidar problem is of course linear. Further, KL is a convex distance particularly appropriate to describe the data discrepancy in the case of Poisson noise. Although the Poisson distribution is a reasonable approximation for describing the lidar experimental observations (such an approximation is particularly reliable at intermediate altitudes while significant systematic corrections occur at low and high altitudes), the input data actually utilized in this paper is the result of a non-linear processing of the experimental measurements and therefore are not Poisson. One of the main results of our investigation is that EM is nonetheless very effective in reconstructing the optical coefficients from such data.

We finally observe that, due to the presence of noise, the maximum likelihood approach tends to provide excessively oscillating solutions, as the random fluctuations of noise are amplified in the inversion procedure \cite{resmerita2007expectation}. However, numerical experiments heuristically show that stopping the iterative process at some optimal iteration regularizes the solution, thus preventing the occurrence of artifacts and over-fitting \cite{benvenuto2014regularization}. Being inspired by results obtained in solar spectroscopy using hard X-ray data \cite{piana2003regularized}, here we use a stopping criterion imposing that the cumulative normalized data residuals satisfy an inequality derived by the central limit theorem.

The plan of the paper is as follows. Section 2 sets up the inverse problem considered and describes EM together with its statistical stopping rule. Section 3 studies the effective resolution associated with the inversion method in the case of noise-free synthetic data. Section 4 performs an error analysis on the profiles used by the European Aerosol Research LIdar NETwork (EARLINET) to assess the inversion performances of different numerical procedures and as training cases for future studies \cite{bosenberg2003earlinet, matthais2004aerosol, bockmann2004aerosol, pappalardo2004aerosol}. Section 5 presents an application to an experimental dataset. Our conclusions are offered in Section 6.

\section{The lidar inverse problem}\label{sec:sec2}

In the atmospheric applications of the lidar technique, measurements of the backscattered mean power $P(z)$ at altitude $z$ and at a specific wavelength $\lambda$ are related to the extinction coefficient $\alpha_{\lambda}(z)$ and the backscatter coefficient $\beta_{\lambda}(z)$ via the lidar signal equation
\begin{equation}\label{a1_E}
P_{\lambda}(z) = \frac{C_{\lambda}}{z^2} \beta_{\lambda}(z) \exp\left(-2 \int_0^z \alpha_{\lambda}(z^{\prime}) dz^{\prime} \right)~,
\end{equation}
where $C_{\lambda}$ is a constant accounting for the laser optical power and the overall detection efficiency, and the laser beam has been assumed to be within the field of view of the receiver at any  altitude (unitary overlap function).
A similar equation can be written for the vibrational Raman lidar signal at wavelength $\mu$ 
\begin{equation}\label{a1}
P_{\mu}(z) = \frac{C_{\mu}}{z^2} \rho(z) \exp\left(- \int_0^z \left( \alpha_{\mu}(z^{\prime}) + \alpha_{\lambda}(z^{\prime}) \right)dz^{\prime} \right)~,
\end{equation}
where $C_{\mu}$ includes laser optical power, detection efficiency and $N_2$ Raman scattering cross section, and the overlap function has also been assumed to be unitary; $\rho(z)$ is the molecular density of the atmosphere, which is assumed to be known, and $\alpha :=\alpha_{\mu}+ \alpha_{\lambda}$ represent the total extinction coefficient, due to atmospheric particles and gases, at $\lambda$ and $\mu$. The molecular contribution can be calculated from the Rayleigh theory of scattering while the relation between the particles' contribution can be derived from Mie scattering theory \cite{ansmann1992independent}. The advantage of Raman equation is evident: since $\rho$ can be determined from a known atmospheric model, $\alpha$ can be computed. Advanced lidar systems with elastic and vibrational Raman channels are available nowadays, and deliver at the same time both backscatter and extinction coefficients. Therefore, in this study we used data coming from a multi-wavelength lidar system (elastic and Raman at several wavelengths). 

In this study we focus on the problem of the inversion of the Raman data; the subscript $\mu$ is omitted but implied.
Simple algebra applied to Eq. (\ref{a1}) leads to 
\begin{equation}\label{Raman_log}
-\log \left( \frac{P(z)\; z^2}{C \;\rho(z)}\right) = \int_{0}^z \alpha(z') dz'~,
\end{equation}
which can be written in compact form as
\begin{equation}
	y = \mathcal{H} x
	\label{eq:ip}
\end{equation}
where:
\begin{itemize}
	\item $\displaystyle y(z)= -\log \left(  \frac{P(z)\; z^2}{C \;\rho(z)} \right)$ is the input data in Eq. (\ref{Raman_log});
	\item $\displaystyle \mathcal{H}: \mathcal{C}\left( [0,z] \right) \to \mathcal{C}\left( [0,z] \right) $ is the integral operator such that $ \displaystyle \left( \mathcal{H}x \right) (z) = \int_0^z x(z')dz'$;
	\item $x$, i.e. the extinction coefficient $\alpha(z)$, is the unknown in Eq. (\ref{Raman_log}).
\end{itemize}

In the following, with a slight abuse of notation, we will be calling $y \in \R^N $ the vector whose $i$-th entry $y_i$ is $y(z_i)$, where $\{z_i~~,~~i=1,\ldots,N\}$ is a set of sampled altitudes where the lidar measurements occur; $x$ is defined accordingly, and $H$ is the $N\times N$ matrix associated with the integral operator $\mathcal{H}$. Therefore we are interested in the numerical solution of the linear system
\begin{equation}\label{linear-system}
y = Hx~.
\end{equation}

The ill-conditioning of Eq. (\ref{linear-system}) implies that straightforward inversion of $H$ would lead to an uncontrollable amplification of noise in the reconstructed extinction profile; therefore regularization is needed to produce stable and physically reliable results \cite{engl1996regularization}. The concept at the basis of regularization consists of replacing the exact problem described by Eq. (\ref{eq:ip}) with the relaxed problem
\begin{equation}
\hat{x} = \arg \min_{x\in \Omega} D(y, x)
\label{eq:ml}
\end{equation}
where $D(y, x)$ is some properly defined measure of the discrepancy between the measured data and the data predicted by the model, and $\Omega$ is a subset (not necessarily proper) encoding the a priori information about the solution. Examples of possible choices for $\Omega$ are: the convex set of all non-negative functions, when one a priori knows that the solution is non-negative; the subspace of functions with compact support, when one a priori knows that the solution is zero outside a closed and limited domain. In constrained maximum likelihood approaches to regularization, the explicit form of the discrepancy measure depends on the statistic of the noise affecting the experimental data. Indeed, maximum likelihood looks for $x$ that maximizes the likelihood, i.e. the probability to observe $y$ from the model $Hx$. In the case of Gaussian noise this probability is
\begin{equation}\label{eq:likelihood-gaussian}
P(y|x) = A \exp(-\|y-Hx\|^2)~,
\end{equation}
where $A$ is a normalization constant and $\|\cdot\|$ is the Euclidean norm. Therefore, in this case, maximizing the likelihood corresponds to minimizing the negative logarithm of Eq. (\ref{eq:likelihood-gaussian}) and $D(x,y)$ becomes the Euclidean discrepancy
\begin{equation}\label{eq:euclidean-discrepancy}
D(y,x) = \|y-Hx\|^2~.
\end{equation}
If one also chooses

\begin{equation}\label{finite-energy} 
\Omega = \{ x \in \R^n ~~|~~\sum_{i=1}^n x_i^2 \leq E \}
\end{equation}
then Eq. (\ref{eq:ml}) and Eq. (\ref{eq:euclidean-discrepancy}) define the Tikhonov regularization method \cite{tikhonov1977solutions}, that can be implemented by solving the minimum problem
\begin{equation}\label{tikhonov}
\hat{x}=\arg\min_{x\in\Omega} \{\|Hx - y\|^2 + \eta \|x\|^2 \}~,
\end{equation}
where $\eta$ is a real, positive regularization parameter to determine via optimization. The Tikhonov method was applied to the first lidar inverse problem in \cite{shcherbakov2007regularized,pornsawad2008ill}. 
Another technique that has been recently applied to regularize the lidar signal equation is described in \cite{pornsawad2012retrieval}; there the authors make use of the Levenberg-Marquardt algorithm, an iterative technique that maximizes a Gaussian likelihood with positivity constraint on the solution; regularization is attained in \cite{pornsawad2012retrieval} by stopping the iterations with the L--curve criterion. 

On the other hand, if the data vector $y$ is a realization of a Poisson random vector, the likelihood can be written as
\begin{equation}\label{likelihood-poisson}
P(y|x) = \prod_{i=1}^N \exp(-(Hx)_i) \frac{(Hx)_{i}^{y_i}}{y_i!}~.
\end{equation}
Again, maximizing this likelihood corresponds to minimizing the discrepancy measure provided by its negative logarithm, known as Kullback-Leibler divergence
\begin{equation}\label{KL}
\displaystyle
D_{KL}(y, x) = \frac{2}{N} \sum_{i=1}^N y_i \log \frac{y_i}{(Hx)_i}  + (Hx)_i - y_i + y_i \log(y_i)~.
\end{equation}
Expectation Maximization (EM) is the constrained maximum likelihood method that chooses $\Omega$ to be $\R^{N+}_*$ (the subset of the strictly positive elements of $\R^N$) and utilizes the Kullback-Leibler divergence as discrepancy measure; it is therefore particularly appropriate in the case of Poisson noise.
The EM algorithm can be derived by using the Karush--Kuhn--Tucker conditions \cite{kuhn1951nonlinear} (named KKT conditions, or theorem) that provide first order necessary conditions for optimality and are largely used for constrained optimization problems \cite{kuhn2014nonlinear}.
In this case, the KKT theorem leads to a fixed point problem that can be turned into the EM iterative algorithm whose standard form is:
\begin{equation}\label{EM-iterative}
	x_{k+1} = \frac{x_k}{H^T 1} \cdot H^T \frac{y}{Hx_k}~~~,
\end{equation}
where $H^T$ is the transpose of $H$ (in the context of EM, $H$ and $H^T$ are often referred to as the forward and the backward operator, respectively); furthermore, with abuse of notation, $1$ is a column vector with all unit entries, and the products and divisions between vectors are intended to be taken component--wise. 

There are several reasons why EM can be specifically appropriate for solving the lidar inverse problem. First, given an initial guess with all strictly positive entries, the solution at each iteration will remain strictly positive, which is a desirable property; we notice that this way to implement the positivity constraint is computationally effective and does not require any projection at each iteration as, for example, in the case of the projected Landweber method \cite{pibe97}. Second, the order of magnitude of the initial guess has no impact on the solution: this can be easily seen by looking at the form of the iteration described by Eq. (\ref{EM-iterative}), where $x_k$ at the right hand side appears both at the numerator and at the denominator. Furthermore, since the problem is a convex optimization problem, in the noise-free case and for an infinite number of iterations, the algorithm will converge to the one and only minimum of the Kullback--Leibler divergence  \cite{lanteri2002penalized}. We finally notice that, as previously observed, EM is based on the use of the Kullback-Leibler divergence which is appropriate in the case of Poisson noise. Lidar measurements obey the Poisson statistic but the data utilized as input in this inverse problem are obtained by processing such measurements (see Eq. (\ref{Raman_log})) and therefore are not exactly Poisson. However, in order to reduce the numerical instability due to ill-conditioning, in the case of noisy data the EM iteration scheme must be stopped according to an optimal criterion, which smoothens the effects of a non-ideal choice of the discrepancy \cite{benvenuto2014regularization}.

The problem of stopping the iterations at a sensible point is known to be a difficult one in the inverse problems literature in general \cite{resmerita2007expectation}, and specifically for the EM method \cite{beetal10,resmerita2007expectation}. In fact,
as the iterations proceed, the EM algorithm produces solutions whose predicted data are closer and closer to the recorded data; at the same time, the solutions become more and more complex, as they attempt to reproduce the finer details in the data; however, as the iterations proceed such fine details are mostly noise. This intuitively explains  why it is so important to stop the iterations within a reasonable range. To achieve that goal we applied a statistical criterion first described in \cite{piana2003regularized}, which stops the iterative scheme when the residual between the measured and the predicted data is compatible with the noise level. More specifically, let
\begin{itemize}
	\item $P_i$ be the measured Raman signal at altitude $z_i$;
	\item $\bar{P}_i$ be the data generated by the extinction profile $x$, i.e. $\displaystyle \bar{P}_i = \frac{C \; \rho}{z_i^2 } \exp\left( -(H x)_i\right)$, at altitude $z_i$;
	\item	$\sigma^P_i$ be the uncertainty of the measurements at altitude $z_i$.
\end{itemize}	
At each iteration, we compute the cumulative mean of the normalized residuals up to $z_i$:
\begin{equation}
\Delta_i = \frac{1}{i} \sum_{j=1}^i \frac{(P_j - \bar{P}_j)}{\sigma_j^P}~~~.
\label{discr_criterion}
\end{equation}
Following \cite{piana2003regularized}, we stop the iterations when $ \displaystyle \Delta_i < \frac{K}{\sqrt{i}} $ for every $i$; specifically we will be using $K=3$.
The rationale behind this criterion is as follows. If the reconstructed solution is close enough to the true value, then the predicted data $\bar{P}$ is close to the exact (noise--free) data. As a consequence, each term in the sum (\ref{discr_criterion}) is the difference between the measured and the exact data, i.e. it is noise. As such, it can be interpreted as an independent variable with zero mean.
By the central limit theorem, as $i$ increases, the distribution of $\Delta_i$ tends to a Gaussian distribution, with variance $\displaystyle \frac{1}{\sqrt{i}}$. It is well known that 99.7\% of the values of a Gaussian fall within 3 standard deviations. In conclusion, if the reconstructed solution is close to the exact value, then the cumulative residuals must fall within the prescribed range; if they don't, then we have not reached a good enough reconstruction, and therefore we proceed with the iterations. 
We stop at the first iteration in which the data generated by the solution is compatible with the measured data, in statistical sense.\\

Before proceeding with the numerical studies, we specify that all the computations described in the present paper have been performed on a MacBook Pro, 3.1 GHz, Intel Core i7; therefore the computational costs below are intended to be relative to this processor.

\section{Effective resolution}

\begin{figure}[!htb]
\begin{tabular}{c}
\includegraphics[width=12cm]{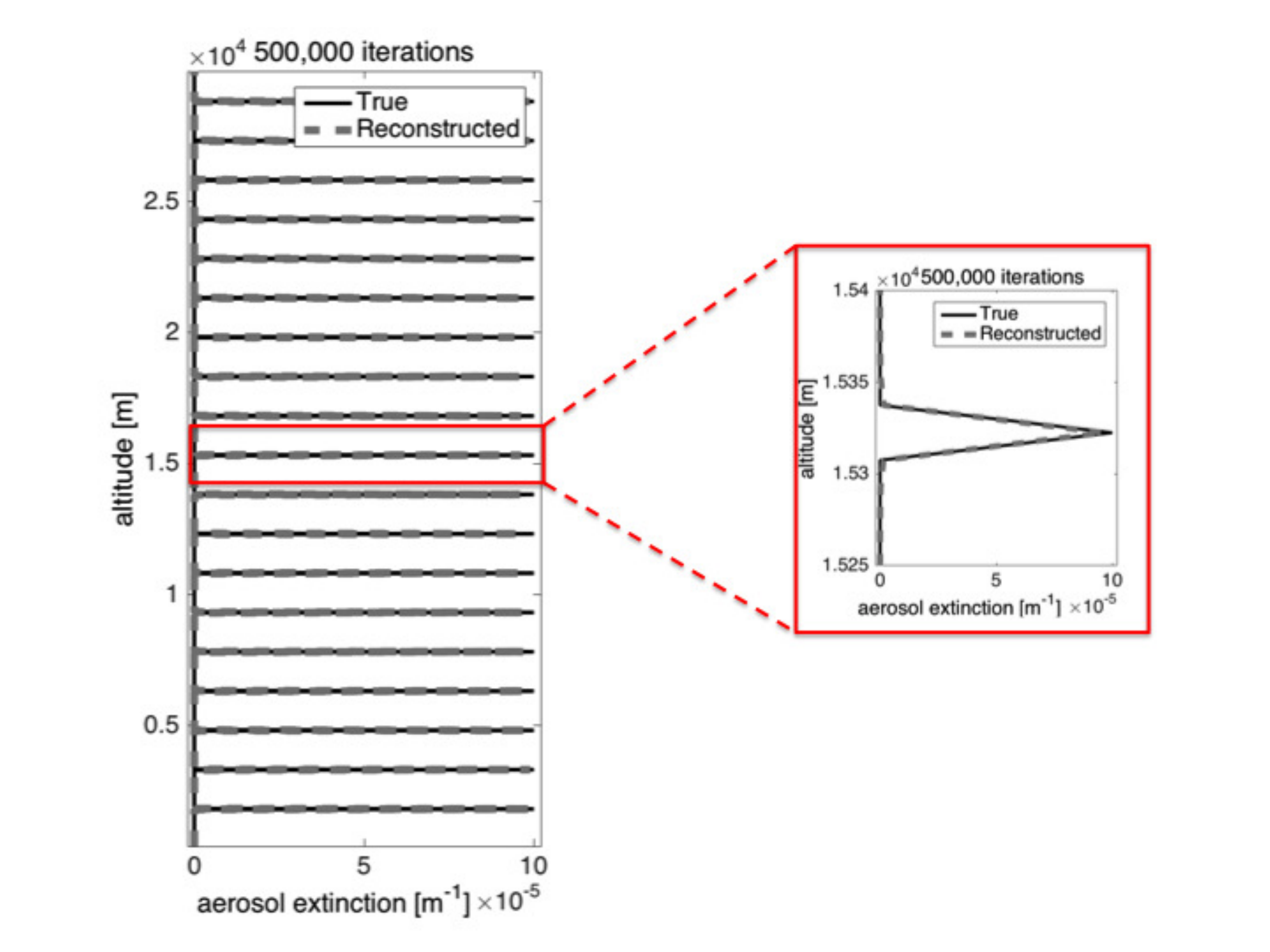} \\
\includegraphics[width=12cm]{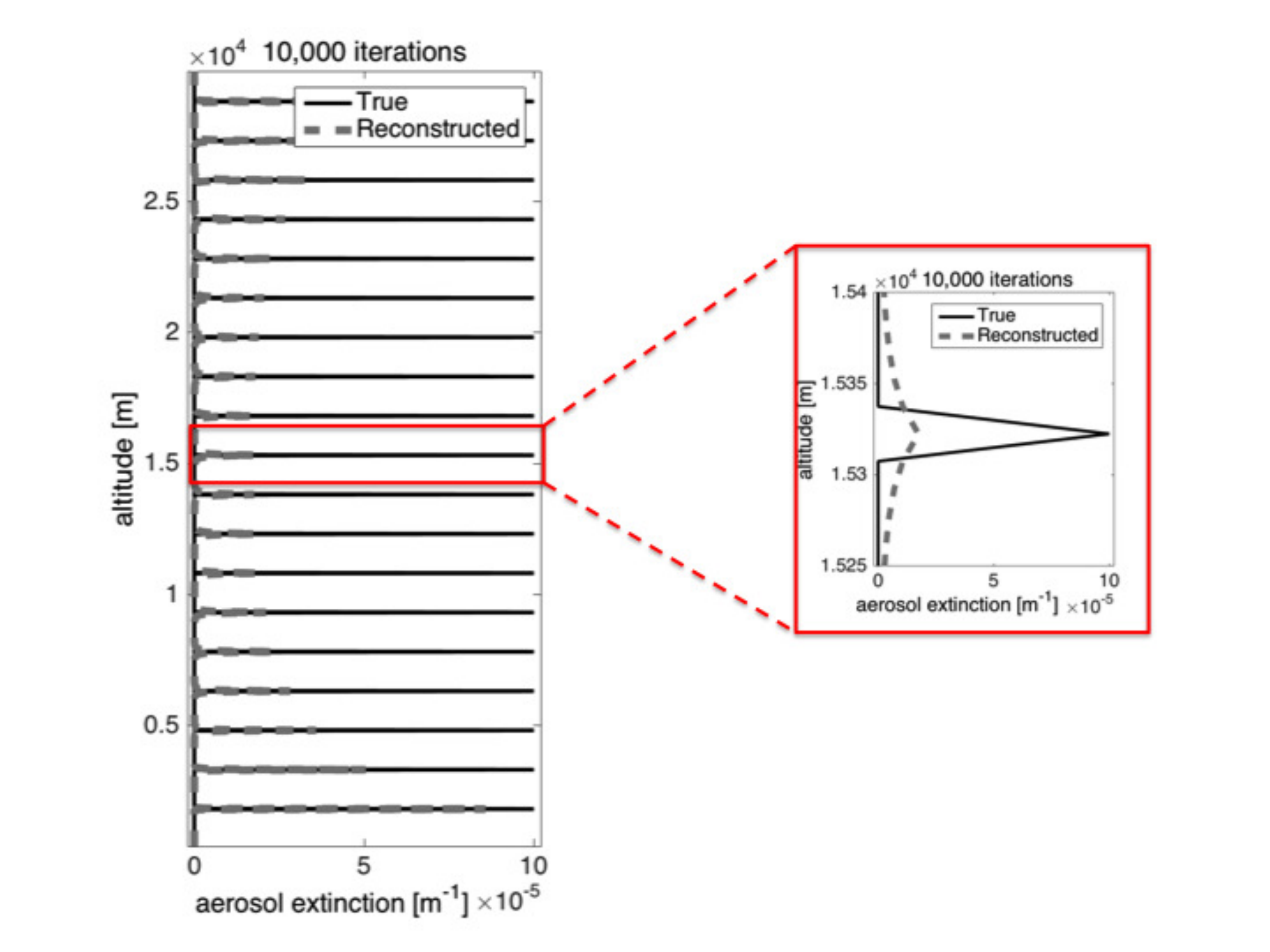} 
\end{tabular}
\caption{Study of the resolution power of EM: reconstruction of a synthetic aerosol extinction profile made of multiple equispaced peaks in the case of 500,000 EM iterations (top panel) and of 10,000 EM iterations (bottom panel). The zoomed plots clearly show that at convergence the resolution power is virtually ideal.}
\label{fig:one}
\end{figure}

\begin{figure}[!htb]
\begin{tabular}{cc}
\includegraphics[width=6cm]{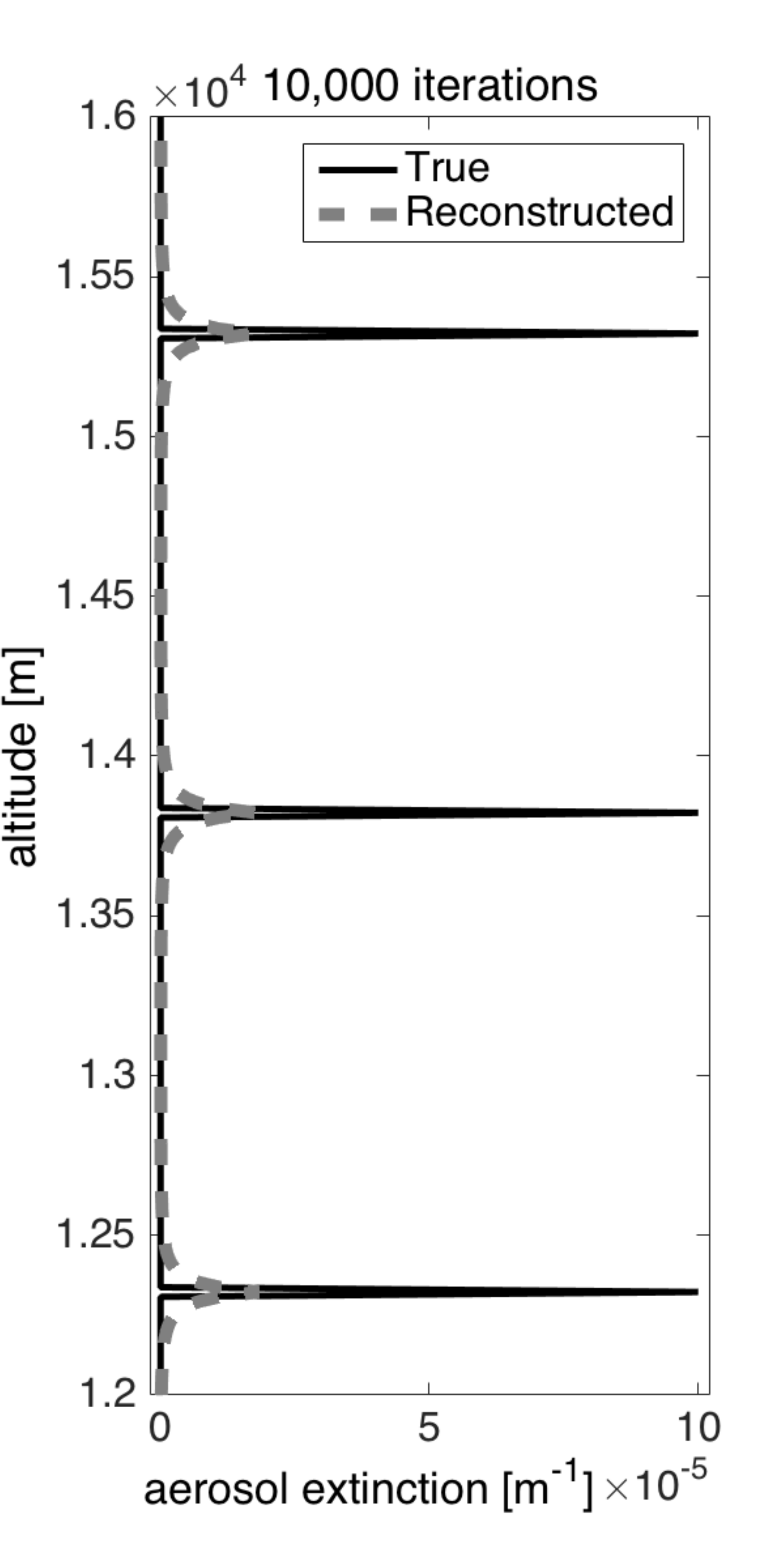} &
\includegraphics[width=6cm]{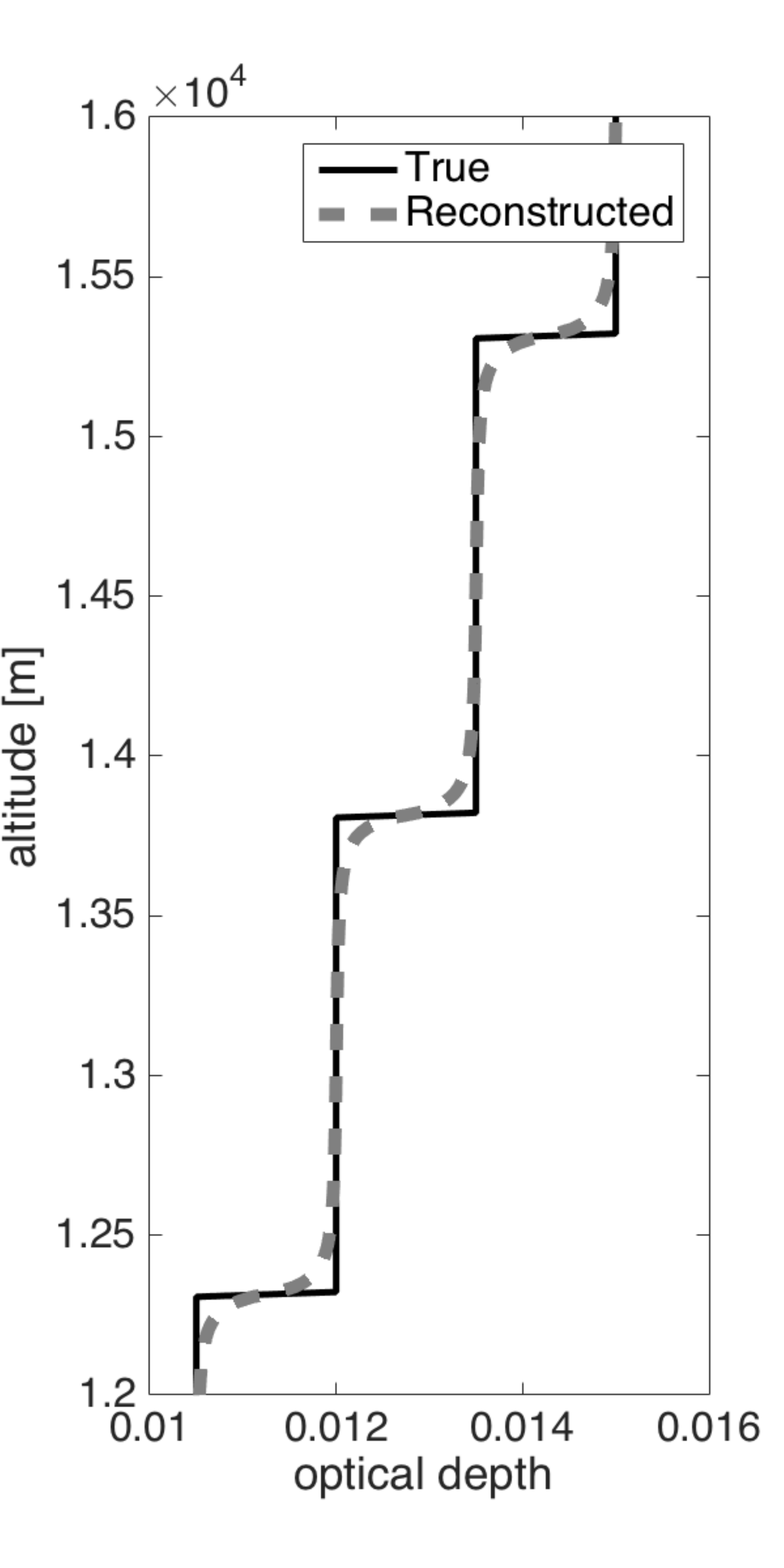} 
\end{tabular}
\caption{Integral (right panel) of a portion of the multiple $\delta$ function profile made of three peaks in the case of $10,000$ EM iterations (left panel).}
\label{fig:two}
\end{figure}

One of the key points in the inversion of lidar data is the quantification of the spatial resolution achievable by means of a specific inversion method.
Indeed, regularization algorithms tend to introduce a blurring effect on the reconstructed backscattering/extinction profiles;
such blurring is, at times, due to the smoothing procedure that is applied directly to the data in order to increase the signal--to--noise ratio;
in other cases, it is due to the smoothness prior used in the penalty term (e.g. in the case of the Tikhonov). Whatever the reason, the result is that
the reconstructed profiles might have significant limitations in distinguishing different peaks when these are closer than a threshold (usually referred to as effective resolution).

Before we proceed with the analysis, we notice that it is not possible to characterize the effective resolution of EM analytically, as one can do, for instance, with the Tikhonov method or Levenberg--Marquardt by means of the averaging kernel procedure. This is due to the non--linear nature of the iterative scheme defined by Eq. (\ref{EM-iterative}). Therefore, in order to investigate the effective spatial resolution of the proposed algorithm, in this section we use synthetic data generated by a simulated
extinction profile containing multiple equispaced peaks (delta functions). 
The simulated profile is discretized every 15 meters from ground level to 15,000 meters. Every ten points (i.e., 150 m) there is a peak of intensity $10^{-4}$ m$^{-1}$; all other points are set to zero. We use this profile to generate synthetic data following Eq. (\ref{linear-system}). %, as shown in Figure \ref{fig:zero} and then apply the EM algorithm.
True and reconstructed profiles are shown in Figure \ref{fig:one} where zooming on a specific peak profile allows visual inspection of the accuracy with which it is reconstructed. We notice that such accuracy is clearly related to the number of iterations utilized by the scheme. Indeed, pushing the algorithm to convergence (for example, letting it work for $500,000$ iterations -- computational time: 37.8 seconds) provides a substantially ideal resolution power. On the other hand, stopping the iterations too early introduces a non-homogeneous blurring that affects the central part of the range more than the extremes. We speculate that the peculiar non--homogeneity of the blurring is due to the nature of the forward and backward operators $H$ and $H^T$, that integrate between zero and $z_k$ and between $z_k$ and $z_{N}$, respectively, and both appear in the EM algorithm. Their effect is therefore more intense in the middle of the range, and symmetric with respect to the central altitude. We also observe that, while failing to reconstruct the height of the peaks, the intensity of the peaks is actually spread over a larger range of altitudes and the overall energy of the signal is conserved. This effect is represented in Figure \ref{fig:two}, which shows the result of integrating three peaks in the reconstruction corresponding to $10,000$ iterations (computational burden of 0.7 seconds). 

We performed many tests, changing the number of peaks and the distance between them. The general behaviour remains the same, with a non--uniform blurring across altitudes, and
better and better reconstructions for more iterations. However, we observed that the strength of the blurring depends on these parameters quite substantially: for example, with three peaks with an inter--distance of 45 meters (3 sampling points) the reconstructed profile is close to the true one already at iteration 20,000; with two peaks with an inter--distance of 150 meters (10 sampling points), 10,000 iterations are enough to obtain almost perfect reconstruction.

At a different level, it is worth to observe that the resolution power of the EM algorithm is also affected by the rearrangements  applied to Eq. (\ref{a1}). 
Such a procedure requires to take the log of $\displaystyle \frac{P(z)z^2}{C\rho(z)}$ and therefore constraints over the values of this quantity are mandatory, to make sure the computation of the logarithm itself makes sense. As a consequence, some of the datapoints are discarded, and this automatically changes the number of altitudes for which it is possible to compute a solution. 

Eventually, our numerical investigations suggest that the problem of quantifying the spatial resolution of EM is a rather complex one, and will be the topic of future work.

\section{Error analysis}

\begin{figure}[!htb]
\begin{center}
	\includegraphics[width=5cm]{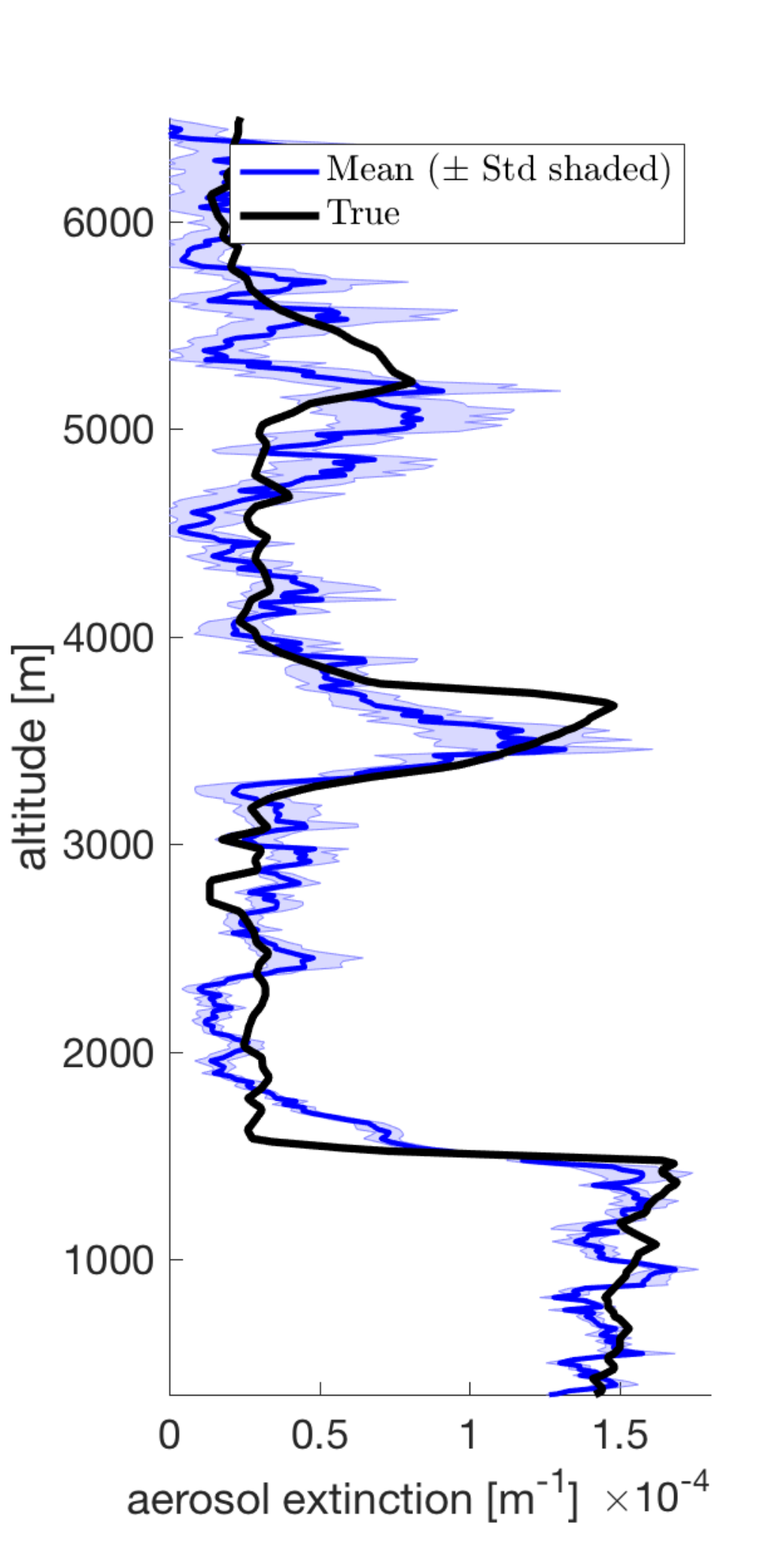}		
	\includegraphics[width=5cm]{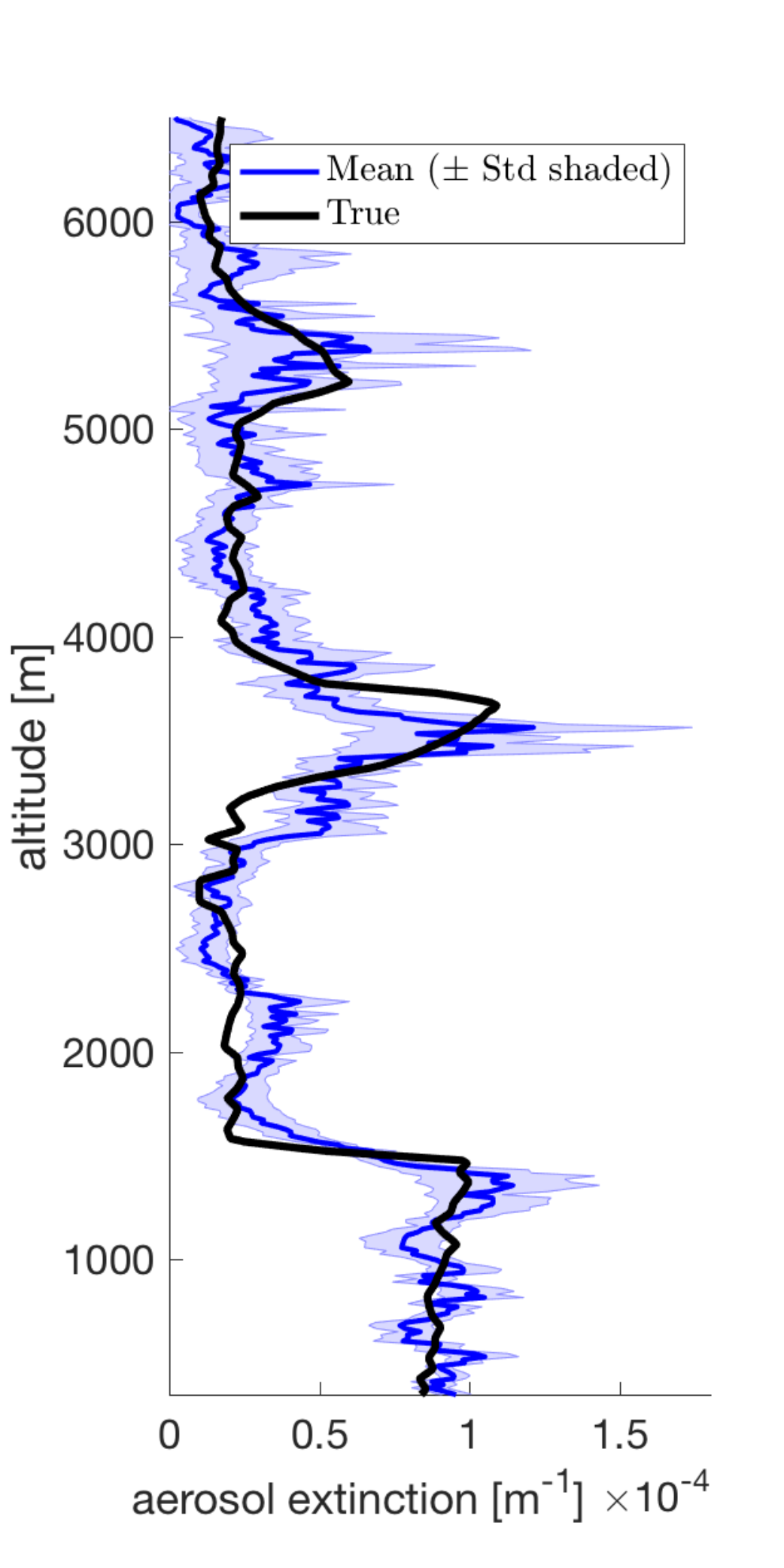}		
	\caption{Mean reconstructed aerosol extinction profile (blue), together with the confidence band and the true profile (black) obtained from the EARLINET synthetic Raman signal, with $K=3$ in the stopping criterion. Left: $\lambda=355$; right: $\lambda=532$.}
	\label{fig:montecarlo}
\end{center}
\end{figure}

We now investigate the impact of noise on the stability of the reconstruction. To do so, we utilize the synthetic extinction profiles used in EARLINET for assessing the performances of different inversion procedures. We consider both wavelengths $\lambda=355$ and $\lambda=532$, and restrict
the range of altitudes to the first 9 kilometers. For each $\lambda$ we construct thirty different realizations of the data, by perturbing with Poisson noise the exact data obtained through Eq. (\ref{linear-system}); then we average these thirty realizations, mimicking the process of recording lidar data every minute, and then averaging over thirty minutes, as typical in EARLINET. We apply EM and obtain an estimated extinction profile. This whole procedure is repeated thirty times to obtain a Monte Carlo variance of the estimated profile. In Figure \ref{fig:montecarlo} we show the mean reconstructed profile (over the thirty repetitions), together with the confidence band and the true profile. In order to stop the iteration we applied the criterion based on the analysis of the cumulative residuals described in Section 2, using $K=3$ (Figure \ref{fig:montecarlo}).  The computational burden, for the entire Monte-Carlo process, is of 21.1 seconds and 19.2 seconds for the $\lambda=355$ and $\lambda=532$ cases, respectively.

\begin{figure}[!htb]
\begin{center}
	\includegraphics[width=5cm]{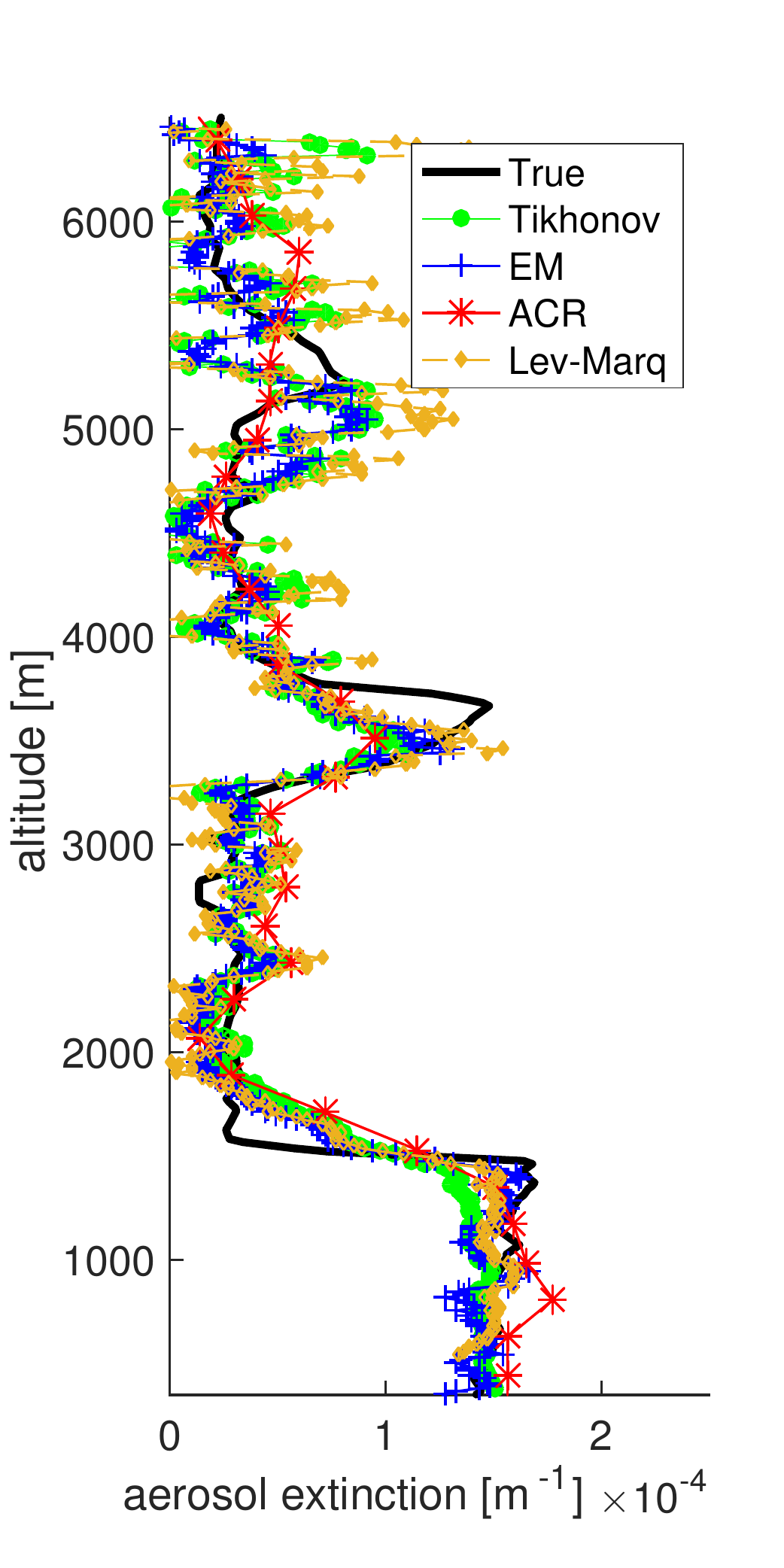}		
	\caption{Comparison of performances between EM, Tikhonov, Levenberg-Marquardt and ACR in the case of the EARLINET data with $\lambda = 355$; for EM, Levenberg-Marquardt and Tikhonov, regularization is realized by exploiting semi-convergence. }
	\label{fig:three}
\end{center}
\end{figure}

As a final test, we compared the performances of EM, of the Tikhonov method, of the Levenberg--Marquardt algorithm \cite{pornsawad2012retrieval} and of a computation tool based on numerical derivative: the ACR (A Complete Retrieval tool for lidar signal analysis) software package from Advanced Lidar Applications srl (\url{http://www.alasystems.it/}).
ACR allows the reconstruction of the extinction profiles by computing the numerical derivative of the Raman signal, making use of several different types of filters.

In this application, both the optimal iteration number for EM and for Levenberg-Marquardt, and the optimal value of the Tikhonov regularization parameter have been selected by means of the criterion of residuals, explained above; for Tikhonov, this is realized by reducing the regularization parameter (rather than proceeding with the iterations) until the criterion is satisfied. On the other hand, the ACR set up requires many parameters to be tuned \textit{ad hoc}, by an experienced user, for optimization between filtering and sharp structure localization according to the signal--to--noise ratio. The included features are:
\begin{itemize}
\item spike rejection: for each set of three consecutive points, the pairwise differences are computed: if they have opposite signs (and bigger that a preselected threshold) the central point is replaced with the average of two adjacent points;
\item binning: the signal can be subjected to a binning of a variable (preselected) number of data points;
\item Gaussian filtering: the signal is subjected to a Gaussian filter of variable (preselected) variance;
\item running average: average current of a variable number of points.
\end{itemize}

Figure \ref{fig:three} shows that EM provides a more stable reconstruction and is more accurate in reproducing the signal peaks with respect to the one obtained via Tikhonov; also the solution provided by Levenberg--Marquardt exhibits substantially more oscillations at higher altitudes. In the figure, for each curve, each point corresponds to an altitude at which the signal has been computed. It is possible to observe that EM outperforms ACR in terms of spatial resolution of the reconstructed data.

\section{Application to experimental data}\label{sec:sec4}

We have used EM to invert two experimental lidar profiles and compared the results with the ones provided by the application of the Tikhonov method, of the Levenberg-Marquardt method, and the application of the numerical derivative method by ACR. 
For this comparison the lidar signals of two distinct apparatuses used in very different atmospheric conditions are considered. The first case (Figure 5) refers to measurements made in Naples with the MALIA (Multiwavelength Aerosol LIdar Apparatus) lidar \cite{boselli2009atmospheric} at 355 nm, with a laser pulse energy of 100 mJ, 20 Hz repetition rate and a 30 cm aperture telescope. Data are collected in Naples 2010/04/20 during a volcanic ash event produced by the eruption of Icelandic volcano Eyjafjallaj\"{o}kull. The second case (Figure 6), is related to measurements made in Naples 2015/04/24 with the AMPLE (Aerosol Multiwavelength Polarization Lidar Experiment) lidar \cite{wang2015calibration} at 355 nm during a Saharan dust event.  The AMPLE lidar makes use of laser source operating at 1000 kHz repetition rate, 0.6 mJ pulse energy, and a receiving telescope of 25 cm diameter. 

\begin{figure}[!htb]
\begin{center}
	\includegraphics[width=4.5cm]{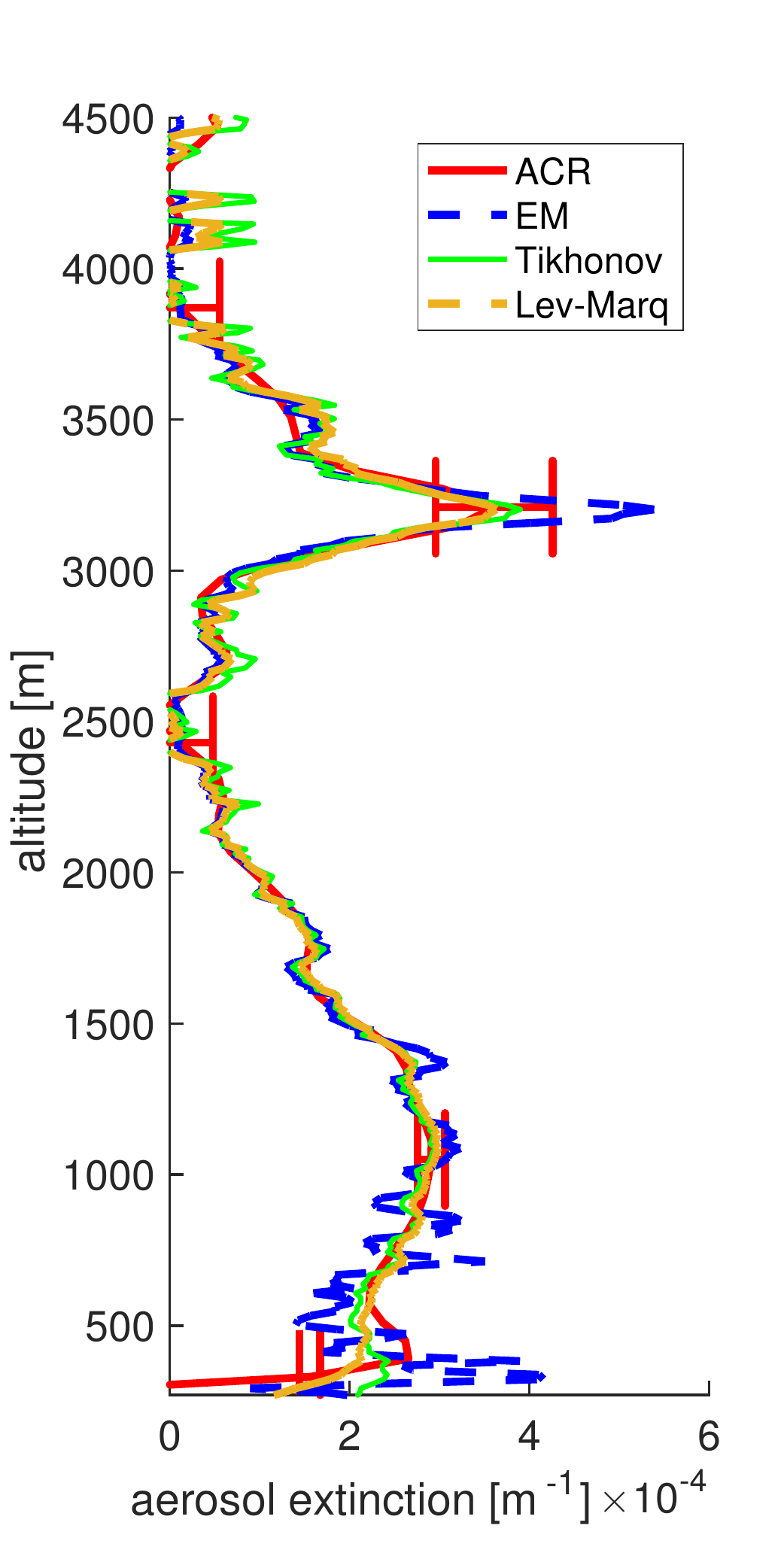}	
	\includegraphics[width=4.5cm]{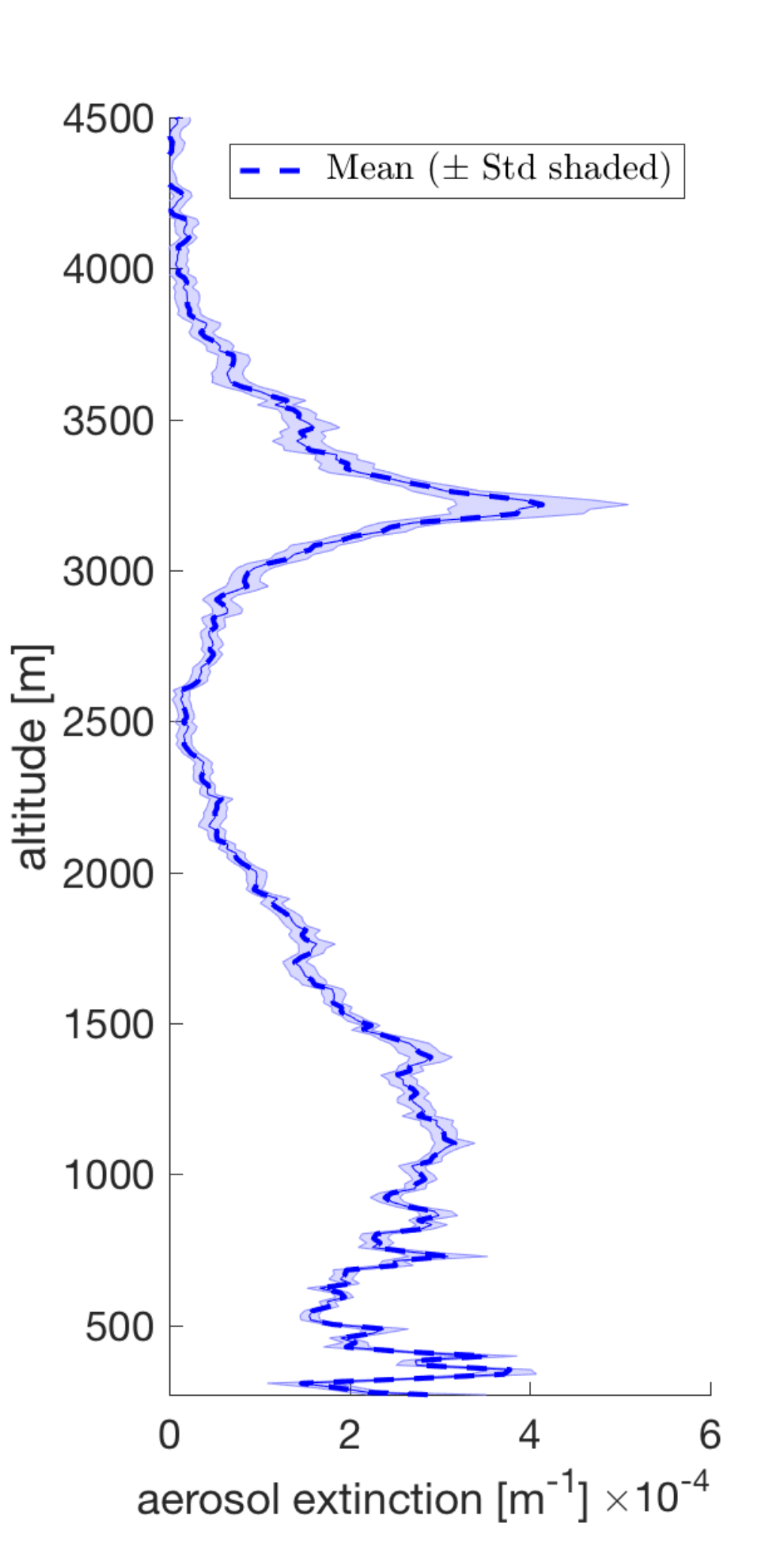}\\
	\includegraphics[width=4cm, trim= 0 -1.3cm 0 0]{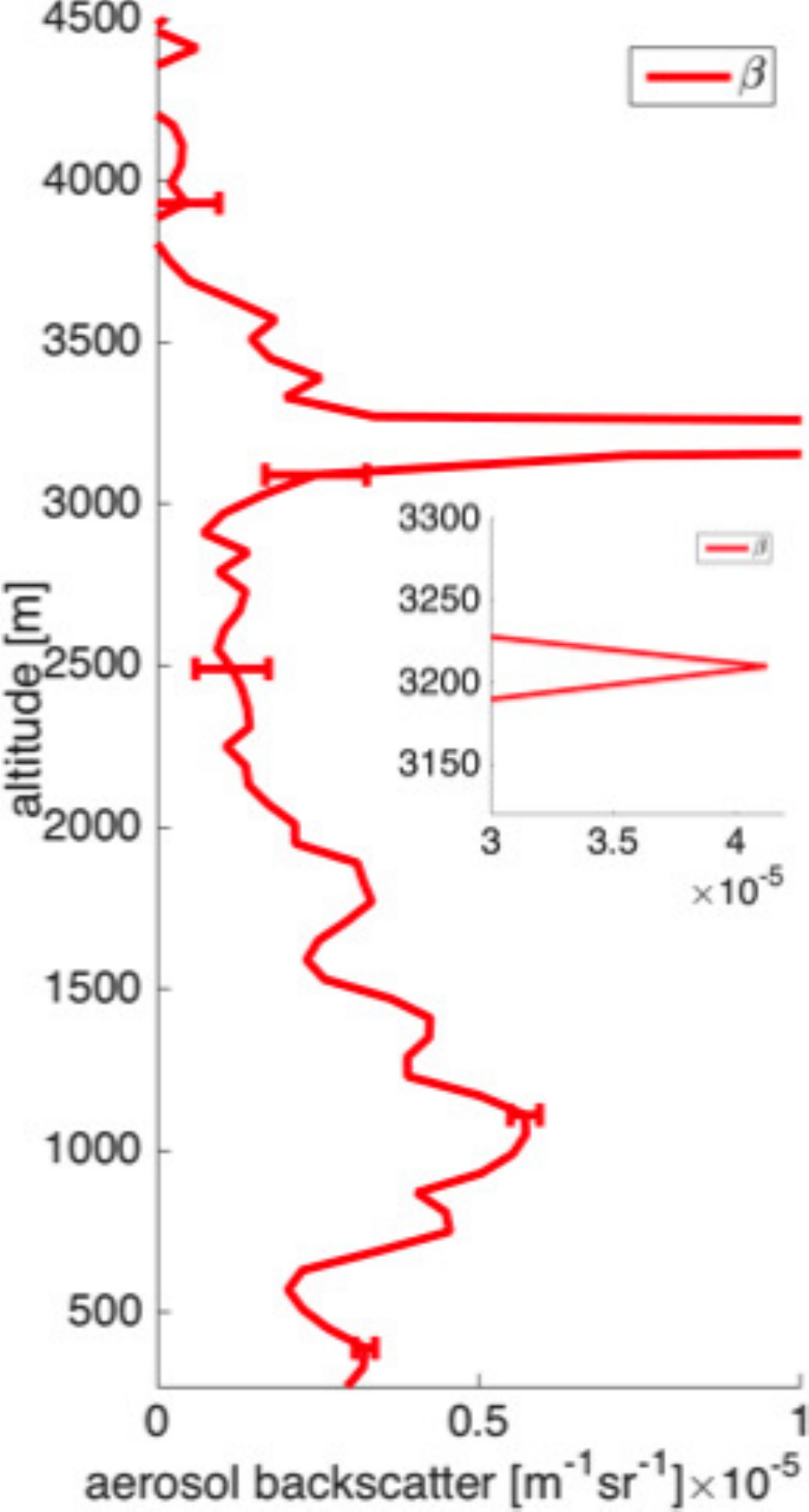}%, 
	\includegraphics[width=4.5cm]{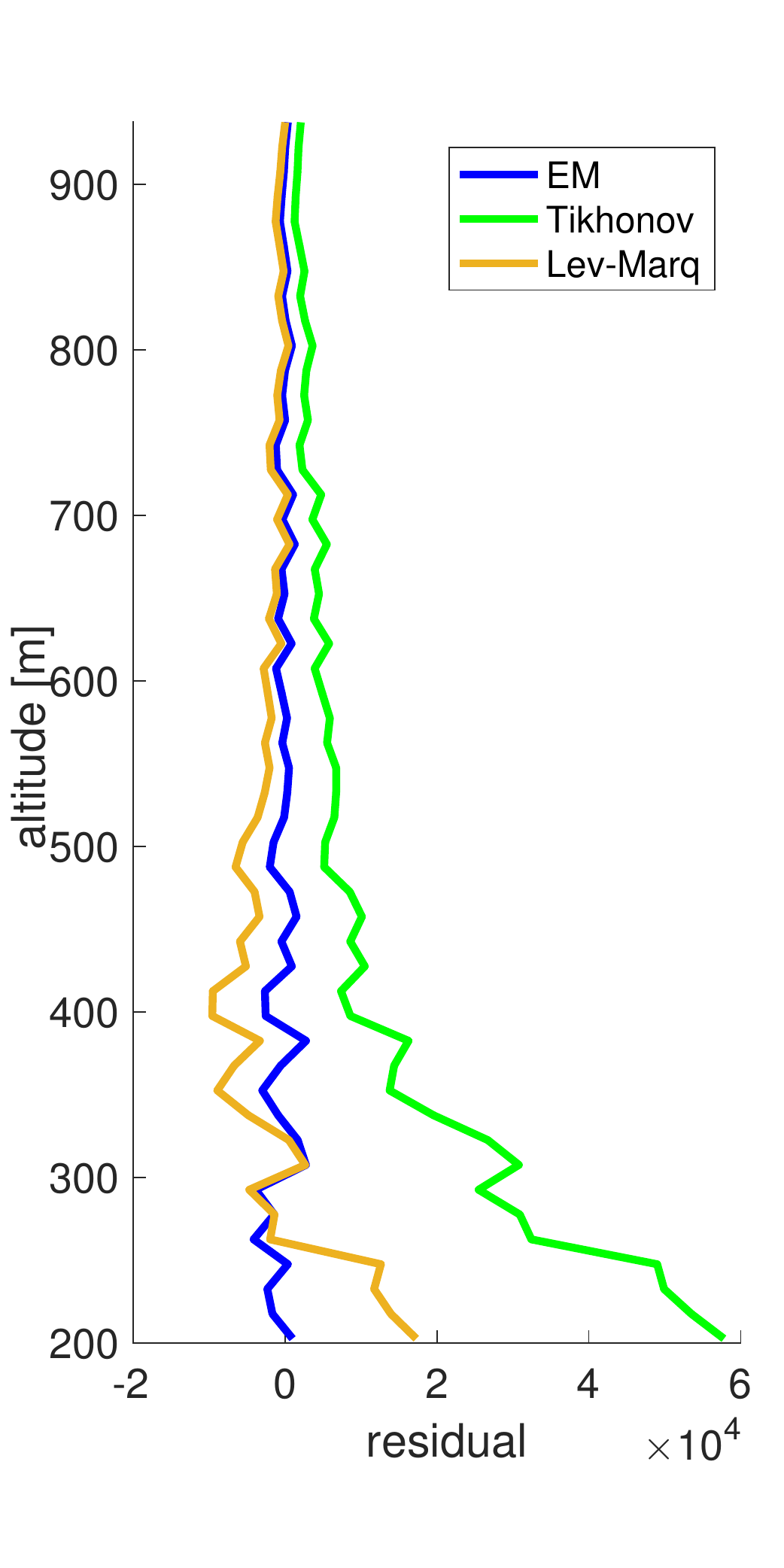}	
	\caption{Top left: analysis of experimental data. The uncertainties on the ACR reconstruction are shown just in some significative points. Top right: EM reconstructed profile with confidence band. Bottom left: backscatter coefficient $\beta$, and a zoom of its peak. Bottom right: residuals (difference between measured and reconstructed signals).}
	\label{fig:four}
\end{center}
\end{figure}

Figure 5 shows, in the top left panel, a comparison of the aerosol extinction profiles obtained by EM, Tikhonov, Levenberg-Marquardt, and numerical derivative by ACR.
In EM, Tikhonov and Levenberg-Marquardt cases, regularization is realized by using the optimality criterion based on the analysis of cumulative residuals. The EM retrieved aerosol extinction profile together with the confidence band is shown in the top right part of the figure. The computational burden for EM, in this case, is of 0.6 seconds, and the selected iteration is 4,601. The Monte--Carlo process employed to obtain the confidence band has a computational burden of 40.9 seconds.

\begin{figure}[htb]
\begin{center}
	\includegraphics[width=6cm]{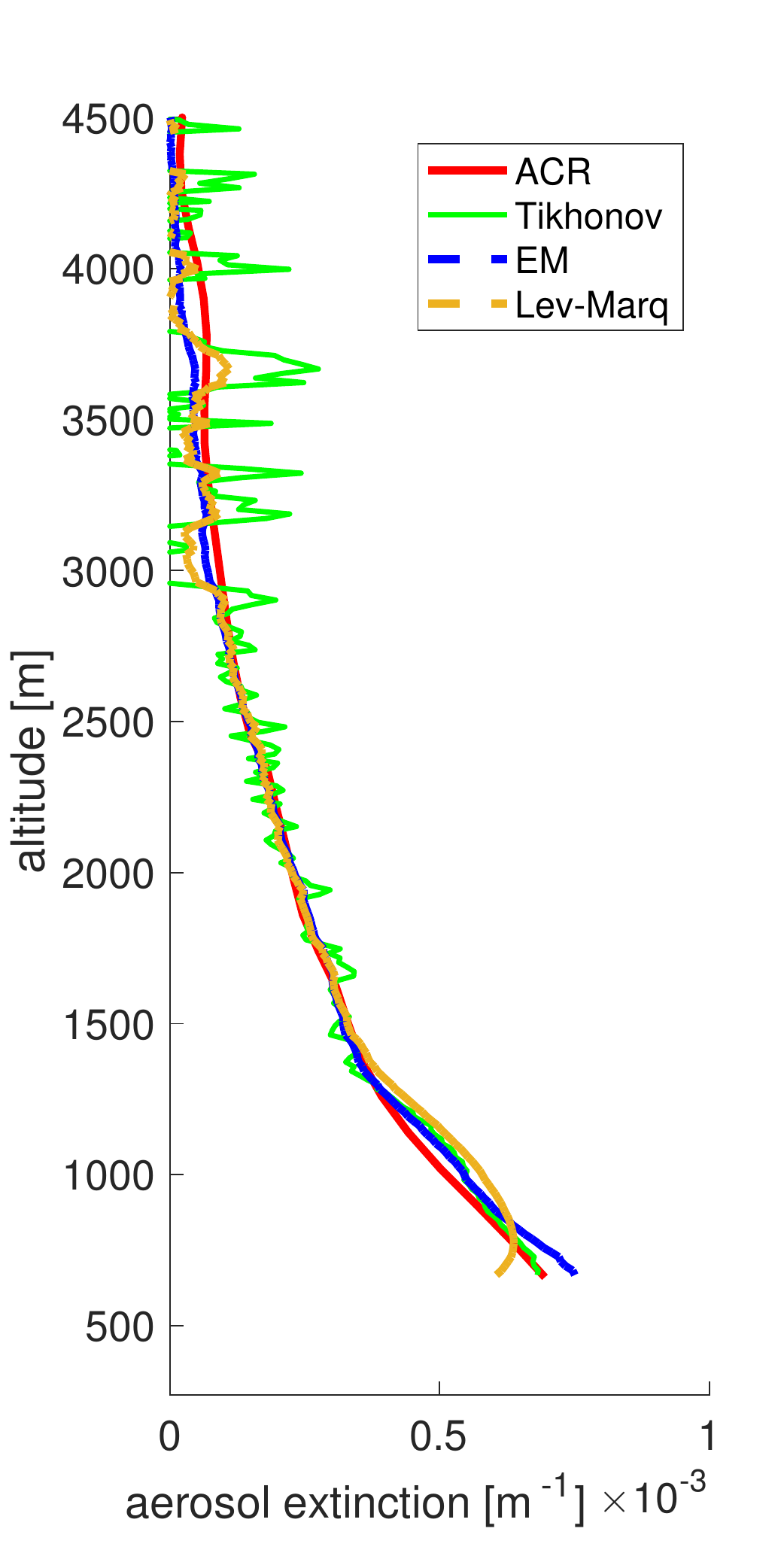}	
	\includegraphics[width=6cm]{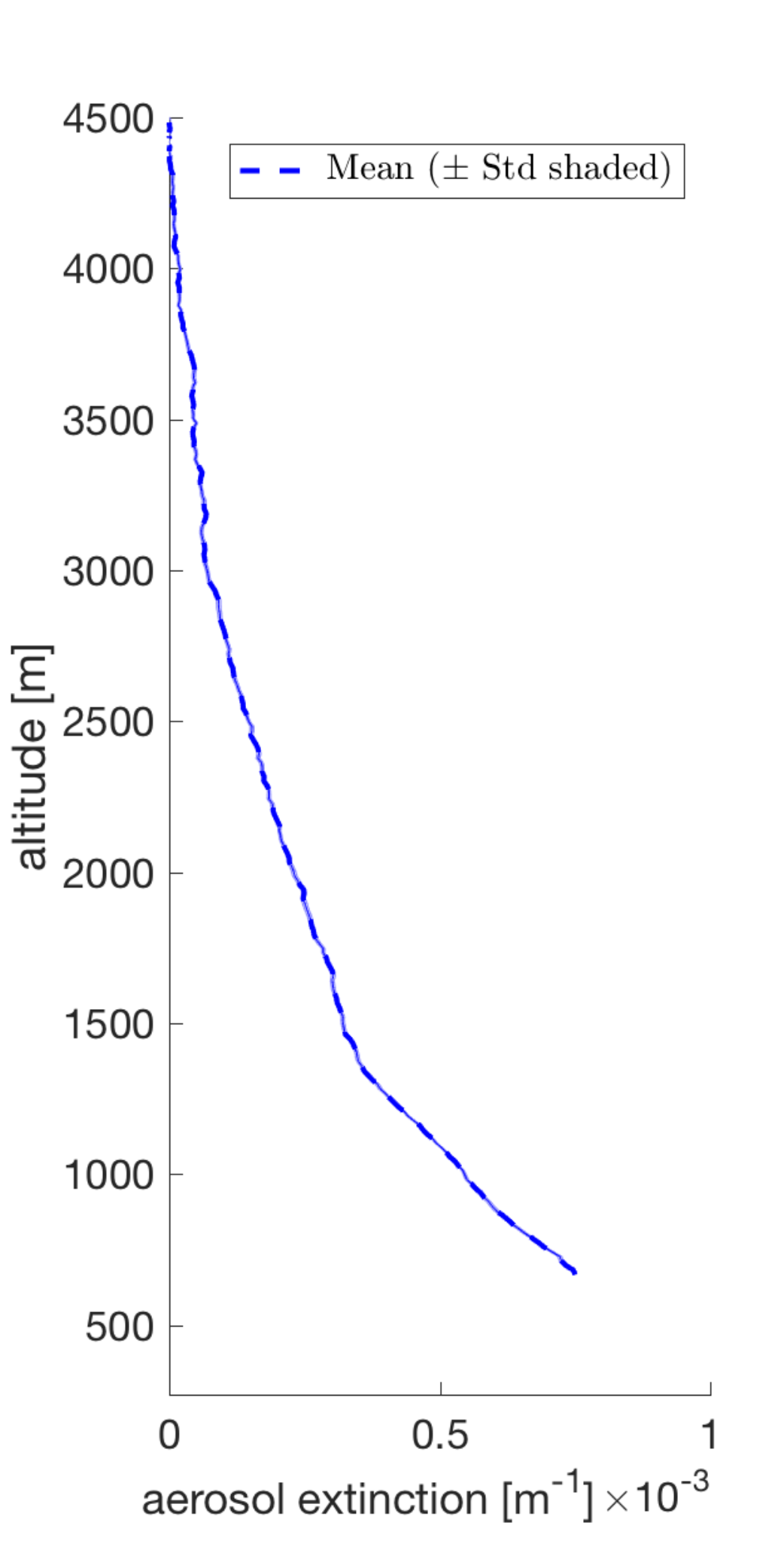}
	\caption{Left: analysis of experimental data. The uncertainties on the ACR reconstruction are shown just in some significative points. Right: EM reconstructed profile with confidence band (very thin in this case).}
	\label{fig:five}
\end{center}
\end{figure}

The Figure points out that EM outperforms both the Tikhonov and the Levenberg-Marquardt method as far as numerical stability is concerned. This is particularly true at higher altitudes, where the signal-to-noise ratio of the experimental measurements notably deteriorates. The EM and ACR solutions are in good agreement each other, taking into account the uncertainty bars on the ACR solution and the confidence bar for EM solution reported in the top right panel. The corresponding aerosol backscattering coefficient vertical profile as retrieved with the Raman method \cite{ansmann1992combined} is plotted in the bottom left panel. This part of the figure clearly shows the presence of a cloud around 3.2 km of altitude. Taking into account that the spatial resolution of the backscattering profile is 60m, the figure demonstrates the superior performance of EM with respect to ACR. In fact, the spatial resolution of the aerosol extinction profile retrieval results to be $180$ m and $300$ m for EM and ACR respectively. The better spatial resolution could be of extreme relevance for the accurate evaluation of the aerosol extinction to backscattering ratio (i.e. the lidar ratio) of the narrow atmospheric features. In addition, we notice that the extinction profile reconstructed by EM presents sharp variations at low altitudes that do not appear in the profiles reconstructed by the other methods. We tend to believe that the reconstruction provided by EM is more reliable than the others, for two main reasons. The first one follows from a visual comparison (Fig. 5, bottom left panel) with the estimated backscatter coefficient (obtained independently from the absorption coefficient, as mentioned above); indeed, the bumps in the backscatter coefficient below 1,500 m are in correspondence with the highest peaks in the extinction coefficient reconstructed by EM. The second reason is that the residuals, i.e. the difference between the measured and predicted data, of the EM solution are smaller and less correlated than those of the other two, as shown in Fig. 5 (bottom right panel), thus suggesting that the EM solution explains the data much better than the others. In this case the comparison is limited to Tikhonov and Levenberg-Marquardt because the solution obtained by ACR has a significantly lower spatial resolution, and comparing the residuals would not be meaningful.

Profiles in Figure 6 have been obtained from measurements performed during a Saharan dust event. The extinction profiles show the presence of aerosol in the low range up to two kilometers altitude. In EM, Tikhonov and Levenberg-Marquardt cases, regularization is realized by using the optimality criterion based on the analysis of cumulative residuals.
Also in this case there is a good agreement between EM and ACR retrieval while instability of Tikhonov results appear above 3 km of altitude. Levenberg-Marquardt appears to be more stable than Tikhonov but still with some fluctuations at higher altitudes. The computational burden for EM, in this case, is of 0.4 seconds, and the selected iteration is 507. The Monte--Carlo process employed to obtain the confidence band has a computational burden of 31.5 seconds.

\section{Conclusions}\label{sec:sec5}

In this paper we have described the application of the EM iterative algorithm, stopped with a cumulative residuals criterion, for reconstruction of the extinction coefficient from Raman lidar data.

Our choice of EM was mostly due to the presence of some desirable properties, namely: the natural and seemless incorporation of a positivity constraint; the fact that, despite being an iterative algorithm, the solution does not depend on the order of magnitude of the first guess; its simplicity and computational efficiency.

We have used synthetically generated data to investigate the capability of EM in reconstructing atmospheric profiles, and observed that the reconstructions provided by EM are more accurate than those provided by the Tikhonov method and by Levenberg--Marquardt. In particular, reconstructions provided by EM seem to have a better overall behaviour across different altitudes: at the lower altitudes they tend to be less over--smoothed than those provided by Tikhonov, where the penalty term encourages continuity of the solution; at higher altitudes they appear to be less affected by noise.

Application of EM to experimental data has been performed considering the lidar signals acquired from two lidar devices differing for laser pulse energy, repetition rate and receiving optics. Further, these two cases correspond to very different atmospheric situations. Both applications confirm that the algorithm can provide meaningful results with real data; in particular, the estimated profiles were similar, overall, to those provided by Tikhonov, by Levenberg--Marquardt and by ACR. But it should be noted that Tikhonov profiles show a greater instability in correspondence of the lowest values of the signal to noise ratio; some oscillations, although fewer than with Tikhonov, are present also in the profile reconstructed by Levenberg--Marquardt. In real cases EM also seems to provide more detailed solutions at low altitudes. To assess the reliability of these additional details we compared the reconstructed profile with an independently obtained profile of the backscatter coefficient, and we also plotted the residuals obtained by the different reconstructions. Both tests indicate that the EM solution is reliable.

In conclusion, both synthetic and experimental data indicate that the solutions provided by EM are comparatively less noisy at high altitudes and more detailed at low altitudes. We speculate that this pleasant feature is mainly due to the discrepancy metric used by EM: indeed, the Kullbach--Leibler divergence, differently from the plain Euclidean norm used by Tikhonov and Levenberg--Marquardt, implies a signal--dependent penalization of the discrepancy between measured and modeled data: stronger signals are allowed to have more variance than small signals. Even though such weighting is optimal for Poisson distributed data (and this is not the case for the log--transformed lidar data), it seems to push the solution to have a better fit at low altitudes while back--projecting less noise at high altitudes. We also speculate that a similar effect might be obtained, in Tikhonov and Levenberg--Marquardt, using a weighted norm in the discrepancy term (corresponding to a Gaussian distribution with non--uniform variance across the altitudes). Future work might be devoted to better investigate this point.

\section*{Funding}\label{sec:sec7}
Advanced Lidar Applications srl.

\section*{Acknowledgments}\label{sec:sec6}
The authors acknowledge EARLINET for providing aerosol synthetic lidar profiles. 

%%%%%%%%%%%%%%%%%%%%%%% References %%%%%%%%%%%%%%%%%%%%%%%%%


\begin{thebibliography}{10}
\newcommand{\enquote}[1]{``#1''}

\bibitem{karyampudi1999validation}
V.~M. Karyampudi, S.~P. Palm, J.~A. Reagen, and H.~Fang,
  \enquote{Validation of the saharan dust plume conceptual model using lidar, meteosat, and ecmwf data,} 
  Bull. Amer. Meteor. Soc. \textbf{80}, 1045 (1999).

\bibitem{di2012raman}
P.~Di~Girolamo, D.~Summa, R.~Bhawar, T.~Di~Iorio, M.~Cacciani, I.~Veselovskii,  O.~Dubovik, and A.~Kolgotin, 
\enquote{Raman lidar observations of a saharan
  dust outbreak event: Characterization of the dust optical properties and determination of particle size and microphysical parameters,} 
  Atmos. Environ. \textbf{50}, 66--78 (2012).

\bibitem{kolev1988lidar}
I.~Kolev, O.~Parvanov, and B.~Kaprielov, 
\enquote{Lidar determination of winds by aerosol inhomogeneities: motion velocity in the planetary boundary layer,}
  Appl. Opt. \textbf{27}, 2524--2531 (1988).

\bibitem{ansmann1990measurement}
A.~Ansmann, M.~Riebesell, and C.~Weitkamp, 
\enquote{Measurement of atmospheric aerosol extinction profiles with a raman lidar,} 
Opt. Lett. \textbf{15}, 746--748 (1990).

\bibitem{pappalardo2004aerosol}
G.~Pappalardo, A.~Amodeo, M.~Pandolfi, U.~Wandinger, A.~Ansmann, J.~B{\"o}senberg, V.~Matthias, V.~Amiridis, F.~De~Tomasi, M.~Frioud, M.~Iarlori, L.~Komguem, A.~Papayannis, F.~Rocadenbosch, and X.~Wang,
\enquote{Aerosol lidar intercomparison in the framework of the earlinet project. 3. Raman lidar algorithm for aerosol extinction, backscatter, and lidar ratio,} 
Appl. Opt. \textbf{43}, 5370--5385 (2004).

\bibitem{pornsawad2012retrieval}
P.~Pornsawad, G.~D'Amico, C.~B{\"o}ckmann, A.~Amodeo, and G.~Pappalardo,
  \enquote{Retrieval of aerosol extinction coefficient profiles from raman lidar data by inversion method,}  
  Appl. Opt. \textbf{51}, 2035--2044 (2012).

\bibitem{povey2014retrieval}
A.~Povey, R.~Grainger, D.~Peters, and J.~Agnew, 
\enquote{Retrieval of aerosol backscatter, extinction, and lidar ratio from raman lidar with optimal estimation,}
 Atmos. Meas. Tech. \textbf{7}, 757--776 (2014).

\bibitem{muller1999microphysical}
D.~M{\"u}ller, U.~Wandinger, and A.~Ansmann, 
\enquote{Microphysical particle parameters from extinction and backscatter lidar data by inversion with regularization: theory,} 
Appl. Opt. \textbf{38}, 2346--2357 (1999).

\bibitem{wang2007retrieval}
Y.~Wang, S.~Fan, and X.~Feng,
 \enquote{Retrieval of the aerosol particle size distribution function by incorporating a priori information,} 
 J. Aerosol Sci.  \textbf{38}, 885--901 (2007).

\bibitem{osterloh2009parallel}
L.~Osterloh, C.~P{\'e}rez, D.~B{\"o}hme, J.~M. Baldasano, C.~B{\"o}ckmann, L.~Schneidenbach, and D.~Vicente, 
\enquote{Parallel software for retrieval of aerosol distribution from lidar data in the framework of earlinet-asos,}
  Comput. Phys. Commun. \textbf{180}, 2095--2102 (2009).

\bibitem{osterloh2013regularized}
L.~Osterloh, C.~B{\"o}ckmann, D.~Nicolae, and A.~Nemuc, 
\enquote{Regularized inversion of microphysical atmospheric particle parameters: Theory and application,} 
  J. Comput. Phys. \textbf{237}, 79--94 (2013).

\bibitem{bertero1998introduction}
M.~Bertero and P.~Boccacci,
\textit{Introduction to Inverse Problems in Imaging} 
(CRC, 1998).

\bibitem{beetal10}
M.~Bertero, P.~Boccacci, G.~Talenti, R.~Zanella, and L.~Zanni, 
\enquote{A discrepancy principle for Poisson data,} 
Inverse Probl. \textbf{26}, 105004 (2010).

\bibitem{resmerita2007expectation}
E.~Resmerita, H.~W. Engl, and A.~N. Iusem, 
\enquote{The expectation-maximization algorithm for ill-posed integral equations: a convergence analysis,} 
Inverse Probl. \textbf{23}, 2575 (2007).

\bibitem{benvenuto2014regularization}
F.~Benvenuto and M.~Piana, 
\enquote{Regularization of multiplicative iterative algorithms with nonnegative constraint,} 
Inverse Probl. \textbf{30}, 035012 (2014).
 

\bibitem{piana2003regularized}
M.~Piana, A.~M. Massone, E.~P. Kontar, A.~G. Emslie, J.~C. Brown, and R.~A. Schwartz, 
\enquote{Regularized electron flux spectra in the 2002 july 23 solar flare,} 
Astrophysical J. Lett. \textbf{595}, L127 (2003).

\bibitem{bosenberg2003earlinet}
J.~B{\"o}senberg and V.~Matthias, 
\enquote{Earlinet: A european aerosol research lidar network to establish an aerosol climatology,} 
Report. Max-Planck-Institut fur Meteorologie \textbf{348}, 1--191 (2003).

\bibitem{matthais2004aerosol}
V.~Matthais, V.~Freudenthaler, A.~Amodeo, I.~Balin, D.~Balis, J.~B{\"o}senberg,  A.~Chaikovsky, G.~Chourdakis, A.~Comeron, A.~Delaval, F.~De Tomasi, R.~Eixmann, A.~Hagard, L.~Komguem, S.~Kreipl, R.~Matthey, V.~Rizi, J.~A. Rodriguez, U.~Wandinger, and X.~Wang,  
\enquote{Aerosol lidar intercomparison in the framework of the earlinet project. 1. Instruments,} 
Appl. Opt. \textbf{43}, 961--976 (2004).

\bibitem{bockmann2004aerosol}
C.~B{\"o}ckmann, U.~Wandinger, A.~Ansmann, J.~B{\"o}senberg, V.~Amiridis,  A.~Boselli, A.~Delaval, F.~De~Tomasi, M.~Frioud, I.~V. Grigorov, A.~Hagard, M.~Horvat, M.~Iarrlori, L.~Komguem, S.~Kreipl, G.~Larcheveque, V.~Matthias, A.~Papayannis, G.~Pappalardo, F.~Rocadenbosch, J.~A. Rodrigues, J.~Schneider, V.~Shcherbakov, and M.~Wiegner, 
\enquote{Aerosol lidar intercomparison in the framework of the earlinet project. 2. Aerosol backscatter algorithms,} 
Appl. Opt.  \textbf{43}, 977--989 (2004).

\bibitem{ansmann1992independent}
A.~Ansmann, U.~Wandinger, M.~Riebesell, C.~Weitkamp, and W.~Michaelis,
\enquote{Independent measurement of extinction and backscatter profiles in cirrus clouds by using a combined raman elastic-backscatter lidar,} 
Appl. Opt. \textbf{31}, 7113--7131 (1992).

\bibitem{engl1996regularization}
H.~W. Engl, M.~Hanke, and A.~Neubauer, 
\textit{Regularization of Inverse Problems} 
(Springer Science \& Business Media, 1996).

\bibitem{tikhonov1977solutions}
A.~N. Tikhonov and V.~I. Arsenin, 
\textit{Solutions of Ill-posed Problems} 
(VH Winston \& Sons, 1977).

\bibitem{shcherbakov2007regularized}
V.~Shcherbakov, 
\enquote{Regularized algorithm for raman lidar data processing,} 
Appl. Opt. \textbf{46}, 4879--4889 (2007).

\bibitem{pornsawad2008ill}
P.~Pornsawad, C.~B{\"o}ckmann, C.~Ritter, and M.~Rafler,
\enquote{Ill-posed retrieval of aerosol extinction coefficient profiles from raman lidar data by regularization,} 
Appl. Opt. \textbf{47}, 1649--1661 (2008).

\bibitem{kuhn1951nonlinear}
H.~W. Kuhn, and A.~W. Tucker,
\enquote{Nonlinear programming,} 
in Proceedings of 2nd Berkeley Symposium. Berkeley: (University of California, 1951), pp. 481--492. 

\bibitem{kuhn2014nonlinear}
H.~W. Kuhn,
\enquote{Nonlinear Programming: A Historical View}
in Traces and Emergence of Nonlinear Programming. (Springer Basel, 2014), pp. 393--414.

\bibitem{shepp1982maximum}
L.~A. Shepp and Y.~Vardi, 
\enquote{Maximum likelihood reconstruction for emission tomography,} 
IEEE T. Med. Imaging \textbf{1}, 113--122 (1982).

\bibitem{pibe97}
M.~Piana and M.~Bertero,
\enquote{Projected Landweber method and preconditioning,} 
Inverse Probl. \textbf{13}, 441--463 (1997).


\bibitem{lanteri2002penalized}
H.~Lanteri, M.~Roche, and C.~Aime, 
\enquote{Penalized maximum likelihood image restoration with positivity constraints: multiplicative algorithms,} 
Inverse Probl. \textbf{18}, 1397 (2002).

\bibitem{boselli2009atmospheric}
A.~Boselli, M.~Armenante, L.~D'Avino, M.~D'Isidoro, G.~Pisani, N.~Spinelli, and X.~Wang, 
\enquote{Atmospheric aerosol characterization over Naples during 2000--2003 earlinet project: Planetary boundary-layer evolution and layering,} 
Bound-Lay Meteorol. \textbf{132}, 151--165 (2009).

\bibitem{wang2015calibration}
X.~Wang, A.~Boselli, A.~Sannino, C.~Song, N.~Spinelli, Y.~Zhao, and C.~Pan,
  \enquote{Calibration of multi-wavelength raman polarization lidar,} 
  EPJ Web Conf. \textbf{89}, 01002 (2015).

\bibitem{ansmann1992combined}
A.~Ansmann, M.~Riebesell, U.~Wandinger, C.~Weitkamp, E.~Voss, W.~Lahmann, and W.~Michaelis, 
\enquote{Combined raman elastic-backscatter lidar for vertical profiling of moisture, aerosol extinction, backscatter, and lidar ratio,}
  Appl. Phys. B-Laser O. \textbf{55}, 18--28 (1992).


\end{thebibliography}
\end{document}